\documentclass[11pt]{amsart}
\usepackage{amsmath,amsthm,amsfonts,amssymb,bm,graphicx}

\usepackage
{hyperref}
\hypersetup{colorlinks=true,citecolor=blue,linkcolor=blue,urlcolor=blue,
pdfstartview=FitH }

\theoremstyle{plain}
\newtheorem{theorem}{Theorem}
\newtheorem{lemma}{Lemma}
\newtheorem{cor}[lemma]{Corollary}

\begin{document}

\author{Valentin Blomer}
\address{Mathematisches Institut, Endenicher Allee 60, 53115 Bonn, Germany} \email{blomer@math.uni-bonn.de}
 
\author{Christopher Lutsko}
\address{Institut f\"ur Mathematik, Universit\"at Z\"urich, Winterthurerstrasse 190, CH-8057 Z\"urich,  Switzerland}\email{christopher.lutsko@uzh.ch}

\newcommand{\Chris}[1]{{\color{red} \sf $\clubsuit\clubsuit\clubsuit$ Chris: [#1]}}
\newcommand{\Valentin}[1]{{\color{blue} \sf $\clubsuit\clubsuit\clubsuit$ Valentin: [#1]}}

  \title{Hyperbolic lattice point counting in unbounded rank}

\thanks{First author supported by DFG through SFB-TRR 358 and EXC-2047/1 - 390685813 and by ERC Advanced Grant 101054336. The second author would like to thank the University of Bonn for hosting him in the summer 2023.
}

\keywords{lattice points, spherical functions, pretrace formula, Eisenstein series}

\begin{abstract} We use spectral analysis to give an asymptotic formula for the number of matrices in ${\rm SL}(n, \Bbb{Z})$ of height at most $T$ with strong error terms, far beyond the previous known, both for small and large rank. 
\end{abstract}

\subjclass[2010]{Primary  11P21, 11N45, 11F72, 22E30}

\setcounter{tocdepth}{2}  \maketitle 

\maketitle
 
\section{Introduction}


\subsection{The main result} Counting lattice points in specified regions dates back at least to Gau{\ss} who gave an asymptotic formula for the number of integer points inside a circle of large radius $R$:
\begin{equation*}
\#\{(x_1, x_2) \in \Bbb{Z}^2 \mid x_1^2 + x_2^2 \leq R^2\} = \pi R^2 + O(R).
\end{equation*}
The error term can be interpreted as the length of the circumference. Following Davenport, this method is now called the Lipschitz principle. By basic harmonic analysis, the error can be improved to $O(R^{2/3})$ (see \cite{Sz}, or \cite[Corollary 4.9]{IK} for a modern treatment).

This question is equally interesting in hyperbolic geometry. The prototypical result goes back to Selberg. Let $\| . \|$ denote the Frobenius norm on ${\rm SL}_n(\Bbb{R})$, i.e.\ $\| g \|^2 = \text{tr}(g^{\top} g)$. Then 
\begin{equation}\label{Se}
\#\{ \gamma \in {\rm SL}_2(\Bbb{Z}) \mid \| \gamma \| \leq T \} = 6 T^2 + O(T^{4/3}).
\end{equation}
While considered classical nowadays, this is nevertheless a difficult result which has never been improved. Selberg's proof was never published. A modern version can be found in \cite[Section 12]{Iw}\footnote{whose proof is written for ${\rm PSL}_2(\Bbb{Z})$ and hence differs by a factor $2$ in the main term}, where the highly non-trivial estimation of the spherical transform is left as an exercise. A long and very different proof of a more general (and only marginally weaker) result can be found in \cite{LP}. 

The aim of this paper is a hyperbolic lattice point count in arbitrary rank.  For $z, w \in {\rm SL}_n(\Bbb{R})$ let $$\mathcal{N}_n(T; z, w) := \#\{\gamma \in {\rm SL}_n(\Bbb{Z}) \mid \| z^{-1} \gamma w \| \leq T\}.$$ This counts  the number of lattice points in the orbit ${\rm SL}_n(\Bbb{Z})w$ in a ball of radius $T$ about $z$. 
Define 
\begin{equation}\label{en}
c_n = \frac{\pi^{n^2/2}}{\Gamma(\frac{n^2 - n + 2}{2}) \Gamma(\frac{n}{2}) \zeta(2) \cdots \zeta(n)}.
\end{equation}
\begin{theorem}\label{thm1} For  $T \geq 1$, $\varepsilon > 0$, $n \geq 3$  and $z, w \in {\rm SL}_n(\Bbb{R})$  we have
$$\mathcal{N}_n(T; z, w) = c_n T^{n(n-1)} + O_{z, w, \varepsilon, n}(T^{n(n-1) -  \delta_n+ \varepsilon})$$
where $$
\delta_3 = 1, \quad \delta_4 = 6/5, \quad \delta_n  = 1 + \frac{1}{\sqrt{2}} + O\Big(\frac{1}{n}\Big).$$
\end{theorem}
The proof works without substantial modifications also for congruence subgroups of ${\rm SL}_n(\Bbb{Z})$. As we will argue below, the numerical values for $\delta_3$ and $\delta_4$ are most likely not improveable by spectral techniques and should be seen as the appropriate generalization of Selberg's bound \eqref{Se}. 

The first result in this direction was proved by Duke-Rudnick-Sarnak \cite[Theorem 1.10]{DRS} and reads  
\begin{equation*}
\delta_n^{{\rm DRS}} =  \frac{1}{n+1}.
\end{equation*}
This was improved only about 25 year later \cite[Theorem 2]{GNY} to
\begin{equation}\label{gny}
\delta_n^{{\rm GNY}} =  \frac{2(n-1)}{(n+1)(n+ \eta)} , \quad \eta = \begin{cases} 0, & n \text{ even},\\ 1, & n \text{ odd}, \end{cases}
\end{equation}
which is asymptotically roughly twice as good. Neither of these bounds recovers \eqref{Se} for $n=2$.  For comparison with Theorem \ref{thm1}, the corresponding savings in  \eqref{gny} are $$\delta_3^{{\rm GNY}} = 1/4,  \quad \delta_4^{{\rm GNY}} = 3/10, \quad \delta_n^{{\rm GNY}} = O(1/n).$$ 
In contrast, our exponent in Theorem \ref{thm1} is uniformly bounded from below. We can slightly improve the asymptotic performance on average over $z$.
\begin{theorem}\label{thm2} Let $\Omega \subseteq {\rm SL}_n(\Bbb{R})$ be a compact set. For  $T \geq 1$, $\varepsilon > 0$ and $n \geq 3$    we have
$$ \int \mathcal{N}_n(T; z, z) dz = c_n \text{{\rm vol}}(\Omega) T^{n(n-1)} + O_{\Omega, \varepsilon, n}(T^{n(n-1) -  \delta_n+ \varepsilon})$$
where $ \delta_n  = 2+ O(1/n).$
\end{theorem}
As we will see below, except for the $O(1/n)$ term this is probably very hard to improve   by spectral techniques.  

\subsection{The methods} In the situation of the classical Gau{\ss} circle problem,  
the strategy is well-known: after a bit of smoothing one applies the Poisson summation formula. The zero frequency yields the main term, the remaining terms are estimated sharply in absolute value. One then optimizes the smoothing parameter to obtain the desired asymptotic formula. Any improvement requires cancellation between the non-zero frequencies (which is possible to some extent  in this particular situation using exponential sum techniques).

In the hyperbolic case, one proceeds similarly. After a bit of smoothing, one applies harmonic analysis in the form of the pretrace formula. The main difference in the non-commutative set-up is that the spectral side looks very different than the geometric side. A sharp estimate of the spherical transform yields \eqref{Se}. Any improvement would require cancellation in the spectral sum over eigenvalues of the hyperbolic Laplacian which is not available with present methods. Consequently, Selberg's error term has never been improved. Note that a Lipschitz principle is not applicable here since the circumference of a large hyperbolic circle has the same order of magnitude as its area. \\

In order to put Theorem \ref{thm1} into perspective, let us first get a feeling for what we can possibly hope for and what might us prevent from obtaining this. The higher rank case offers two major difficulties: 

a) a sharp estimate for the spherical transform of the characteristic function on a ball or a smoothed version thereof. This is a problem in the analysis on Lie groups where we obtain sharp bounds (at least on the tempered spectrum); cf.\ Section \ref{spher}. 

b) an understanding of the non-tempered spectrum. This is deeper than it might sound. First of all, in higher rank there could be infinitely many linearly independent cusp forms violating the Ramanujan conjecture. More seriously, residual Eisenstein series are known to violate the Ramanujan conjecture, some of them quite drastically, and will therefore contribute a large error term. One could hope that such Eisenstein series make up only a small portion of the spectrum which compensates their degree of non-temperedness. However, an analysis of the pretrace formula also asks for pointwise bound of all appearing spectral components, and in higher rank we have very little control on the sup-norm of automorphic forms, regardless of whether they are cuspidal or Eisenstein. The best we can do in general is to use the pretrace formula backwards and estimate the sup-norm by the square-root of the spectral density. In this way, however, we sacrifice a large portion of our  knowledge on the sparsity of such Eisenstein series. The reader may notice from the proof that on the other hand two features are working in our favor. On the one hand residual Eisenstein series lie on many Weyl chamber walls which reduces the spectral density (somewhat). On the other hand, the spherical transform of the characteristic function of the set of matrices with norm $\leq T$ behaves (slightly) better on the non-tempered spectrum. \\

Let us make the preceding remarks a bit more quantitative. We smooth out the characteristic function on $[0, T]$ to a function with support in $[0, T(1 + \delta)]$ where $\delta$ will be chosen later as a function of $T$. This introduces an error in the main term of $O(T^{n(n-1)} \delta)$. We will show in Section \ref{spher} below that the spherical transform of such a function decays roughly like
$$\frac{T^{\frac{1}{2}n(n-1) + n \| \Re \mu \|} }{\| \mu \|^{\frac{1}{4} n(n+1)}} (1 + \delta \| \mu \|)^{-A}$$
for any $A > 0$ and a spectral parameter $\mu = (\mu_1, \ldots, \mu_n)$;   we follow the convention the tempered spectrum satisfies\footnote{Some author prefer to normalize the tempered spectrum such that $\mu \in \Bbb{R}^n$.} $\mu \in (i\Bbb{R})^n$. Here $\| . \|$ denotes the maximum norm of a vector. This bound is an oversimplification and holds only if $\mu$ is in generic position, i.e.\ away from the Weyl chamber walls, but let us proceed anyway. 
Recall that by Weyl's law \cite{Mu} there are $O( \delta^{1- n(n+1)/2 })$ linearly independent cusp forms with  $\mu \ll \delta^{-1}$, the point after which the spherical transform becomes negligible. 

If we ignore the non-tempered spectrum as well as the problem of pointwise bounds for automorphic functions, we obtain a total error term
\begin{equation}\label{best}
\ll T^{n(n-1)} \delta + T^{\frac{1}{2}n(n-1)} \delta^{1- \frac{1}{2}n(n+1) + \frac{1}{4} n(n+1) } \ll T^{n(n-1) - 2\frac{n-1}{n+1}}
\end{equation}
upon choosing $\delta = T^{-2(n-1)/(n+1)}$. 
This is the most optimistic generalization of Selberg's argument to higher rank and recovers \eqref{Se} for $n=2$. Since even for $n=2$ this has been the state of the art for more than half a century, it seems essentially impossible to go beyond such a bound. We emphasize that unlike \eqref{gny}, this suggests a saving of \emph{asymptotically constant} size.\\


This is precisely what we achieve in Theorem \ref{thm1} with a constant $1 + 1/\sqrt{2} \approx 1.7$.  The cases $n=3$ and $n=4$ match our (rather optimistic) heuristic \eqref{best}.  In Theorem \ref{thm2} we reach asymptotically the ``best-possible'' constant $2$ on average over $z$. (For $n=2$, the paper \cite{PR} shows that a rather sophisticated additional application of the Kuznetsov formula   can exploit the average a bit more strongly, but this seems very hard to implement in higher rank.)

 We remark that our main results are somewhat similar in spirit and quality (relative to previous results) as the recent paper \cite{JK}, however with the major difference that Theorem \ref{thm1} is a pointwise result, valid for all (fixed) $z, w$, whereas \cite{JK} holds for almost all points. 




\subsection{Acknowledgments}
We thank Peter Sarnak for the discussion which led to this project, and Alex Kontorovich for insightful discussions. 

\section{Preliminaries}

\subsection{Some notation}
We write $G = {\rm SL}_n(\Bbb{R})$, $\Gamma = {\rm SL}_n(\Bbb{Z})$, $K = {\rm SO}(n)$, $A \subseteq {\rm SL}_n(\Bbb{R})$ for the group of diagonal matrices with positive entries and determinant 1, $N$ for the group unipotent upper triangular matrices, $W$ for the Weyl group.  We decompose $G = KAN = NAK$. Accordingly, we have the Iwasawa projections\footnote{The double use of $A$ will not lead to confusion.} $A, H : G \rightarrow \mathfrak{a}$ such that $g \in K \exp(H(g)) N= N\exp(A(g)) K$. We denote by $dn, da, dk, dg$ the usual Haar measures, normalized as in \cite[Appendix A]{DRS}\footnote{but the particular normalization plays no role for the purpose of this paper}, in particular $\text{vol}(K) = 1$. As usual, we write $\rho = (\frac{n-1}{2}, \frac{n-3}{2}, \ldots, \frac{3-n}{2}, \frac{1-n}{2}) \in \mathfrak{a}^{\ast}$ and $C_{\rho}$ for the convex hull of the points $\{w \rho \mid w \in W\}$. We often identify $\mathfrak{a}_{\Bbb{C}}^{\ast} \cong \{\mu \in \Bbb{C}^n  \mid \sum \mu_j = 0\}$ and  equip vectors in $ \mathfrak{a}^{\ast}_{\Bbb{C}}$ with the max-norm $\| . \|$. We recall that all relevant spectral parameters $\mu \in \mathfrak{a}^{\ast}_{\Bbb{C}}$ satisfy
\begin{equation}\label{Lambda}
    \sum_{j=1}^n \mu_j = 0, \quad \{\mu_1, \ldots, \mu_n\} = \{-\bar{\mu}_1, \ldots, -\bar{\mu}_n\}, \quad \mu \in i\mathfrak{a}^{\ast} +  C_{\rho}. 
\end{equation}
In the following we regard $n \geq 2$ as fixed, and all implied constants may depend on $n$.

\subsection{Spherical transform and pretrace formula} For $\mu \in \mathfrak{a}^{\ast}_{\Bbb{C}}$ we define the spherical function (cf.\ \cite[p.\ 418 \& p.\ 435]{He} with $\mu$ in place of $i\lambda$)
\begin{equation*}
\phi_{\mu}(g) = \int_K e^{(-\rho + \mu)H(gk)} dk= \int_K e^{(\rho + \mu)A(kg)} dk.
\end{equation*} 
We have 
\begin{equation}\label{trivphi}
\phi_{\rho}(g) = 1, \quad |\phi_{\mu}(g)| \leq 1
\end{equation}
for all $g\in G$ and $\mu$ satisfying \eqref{Lambda}. 
For  a (measurable) compactly supported, bi-$K$-invariant function $f : G \rightarrow \Bbb{C}$ we define the spherical transform (cf.\ \cite[p.\ 449]{He})
$$\tilde{f}(\mu) = \int_G f(g) \phi_{-\mu}(g) dg.$$
It is well-known (\cite[p.\ 450]{He}) that this is the composition of the Abel transform  
$$\mathcal{A}_f(a) = e^{\rho(\log a)} \int_{N} f(an)dn, \quad a \in A,$$
and the Fourier transform 
\begin{equation}\label{abel}
\tilde{f}(\mu) = \int_A \mathcal{A}_f(a) e^{-\mu \log a} da.
\end{equation}
If $f_1, f_2$ are two (measurable) compactly supported, bi-$K$-invariant functions, the convolution 
$$(f_1 \ast f_2)(x) = \int_G f_1(x g^{-1})f_2(g) dg$$
satisfies \cite[p.\ 454]{He} 
\begin{equation}\label{conv}
\widetilde{f_1 \ast f_2} = \tilde{f_1} \tilde{f_2}.
\end{equation}
The spherical transform comes up in the pretrace formula \cite{Sel}: for a smooth, compactly supported, bi-$K$-invariant function $f : G \rightarrow \Bbb{C}$ and $z, w\in \Gamma \backslash G/K$ the spectral expansion of the automorphic kernel on the left hand side (cf.\ \cite[(2.5)]{Sel}) of the following display together with the uniqueness principle (cf.\ \cite[(1.8)]{Sel}) is
\begin{equation}\label{pretrace}
   \sum_{\gamma \in \Gamma} f(z^{-1} \gamma w) = \int  \tilde{f}(\mu_{\varpi}) \varpi(z) \overline{\varpi(w)} d \varpi.
\end{equation}
The notation should be interpreted as follows: the right hand side runs over cusp forms and Eisenstein series (including residual Eisenstein series) for the group $\Gamma$, and $d\varpi$ is the counting measure on the discrete spectrum. Each spectral parameter $\mu_{\varpi} \in  \mathfrak{a}_{\Bbb{C}}^{\ast}$   of $\varpi$ is only defined modulo the action of the Weyl group $W$. We parameterize the spectrum in detail in Subsection \ref{Eisen}.     

\subsection{Sup-norm bounds}
We apply the pretrace formula   in the other direction to obtain the following (generic)   bounds for automorphic forms appearing on the right hand side of \eqref{pretrace}.

\begin{lemma}\label{lem1}
Let $\mu \in i \mathfrak{a}^{\ast}$, $B(\mu)  \subseteq \mathfrak{a}_{\Bbb{C}}^{\ast}$  a ball of size $O(1)$ about $\mu$ and $z \in G$. Then we have
\begin{equation*}
\int_{B(\mu)} |\varpi(z)|^2 d\varpi \ll_{z}  \prod_{1 \leq i < j \leq n}(1 + |\mu_i - \mu_j|).
\end{equation*}
\end{lemma}

\textbf{Proof.} This is well-known and implicit for instance in \cite{BM}. For convenience we recall the proof. Let  $f$ be a fixed   function on  $\mathfrak{a}_{\Bbb{C}}^{\ast}$ with compactly supported Fourier transform such that $f$ is  real  on $i\mathfrak{a}^{\ast}$,   $\Re f $ is non-negative in the strip  $\{\lambda\in \mathfrak{a}_{\Bbb{C}}^{\ast} :  | \Re \lambda_j | \leq  \| \rho\|\}  $,  and $\Re f \geq 1$ on a ball in $\mathfrak{a}_{\Bbb{C}}^{\ast}$ about $0$ of radius $\| \rho \|_2$. Such a function was constructed explicitly in \cite[Lemma 1]{BM}  and the subsequent display. Then we choose
$$\tilde{f}_{\mu}(\lambda) := \Bigl( \sum_{w \in W}f(\mu - w \cdot \lambda)\Bigr)^2.$$
This again has compactly supported Fourier transform, and  the support   is independent of $\mu$. 
One verifies quickly that 
$$\tilde{f}_{\mu}(\lambda) \geq 0$$
for all $\lambda$ satisfying \eqref{Lambda} and $$\tilde{f}_{\mu}(\mu) \geq 1.$$ 
Moreover, the rapid decay along the real axis shows 
$$\tilde{f}_{\mu}(\lambda) \ll_A \max_{w\in W}(1 + \| \Im \mu - w\cdot  \lambda \|)^{-A}$$ for $\lambda \in i\mathfrak{a}^{\ast}$ and any $A > 0$.  By the Harish-Chandra inversion formula featuring the Harish-Chandra $\textbf{c}$-function together with the trivial bound \eqref{trivphi} for elementary spherical functions, we see that the inverse spherical transform $f_{\mu}$ of $\tilde{f}_{\mu}$ has compact support and  satisfies the bound 
\begin{equation*} 
f_{\mu}(g) \ll  \prod_{1 \leq i < j \leq n}(1 + |\mu_i - \mu_j|). 
\end{equation*}
We now conclude from \eqref{pretrace} that
 $$\int_{B(\mu)} |\varpi(z)|^2 d\varpi \ll \int_{B(\mu)} |\varpi(z)|^2 \tilde{f}_{\mu}(\mu_{\varpi}) d\varpi = \sum_{\gamma\in \Gamma} f_{\mu}(z^{-1}\gamma z) \ll  \prod_{1 \leq i < j \leq n}(1 + |\mu_i - \mu_j|)$$
 as desired. \\
 
In Subsection \ref{Eisen} we give a better bound on average over $z$, taken from \cite{JK2}.

\subsection{Smoothing} 
For $T > 1$ we write $\chi_T : K\backslash G/K \rightarrow \Bbb{R}_{\geq 0}$ for the characteristic function on $\| g \| \leq T$. 
For $0 < \delta< 1$ and let $\psi_\delta : K\backslash G/K \rightarrow \Bbb{R}_{\geq 0}$ be a smooth $L^1(AN)$-normalized function  supported in $B_{\delta}  := \{ g \in G \mid \max(\| g \|_2, \| g^{-1}\|_2) \leq 1 + \delta\}$ where $\|. \|_2$ is the matrix norm induced from the Euclidean vector norm. Let
\begin{equation}\label{convdef}
\chi_{T, \delta}  := \chi_T \ast \psi_{\delta}. 
\end{equation} 
By \cite[Lemma 3.3]{DRS} we have for $hg^{-1} \in \text{supp}(\chi_T)$ and $g\in \text{supp}(\psi_{\delta})$ the inequalities
$$T \geq \| h g^{-1} \|  \geq  \frac{\| h g^{-1} g \|}{\| g \|_2} \geq \frac{\| h  \|}{1 + \delta} .$$
On the other hand we have for $\| h \| \leq T(1+\delta)^{-1}$ and $g\in \text{supp}(\psi_{\delta})$  that
$$ \| h g^{-1} \| \leq \| h \| \| g^{-1} \|_2 \leq T.$$
These two inequalities imply $\chi_{T(1 + \delta)^{-1}} \leq \chi_{T, \delta} \leq \chi_{T(1+\delta)}$ and hence
\begin{equation}\label{ineq}
\chi_{T(1+\delta)^{-1}, \delta} \leq \chi_T \leq \chi_{T(1+\delta), \delta}. 
\end{equation}

\subsection{Parametrization of the spectrum}\label{Eisen} We describe the various types of Eisenstein series in classical language. See \cite{MW1, MW2} for a detailed description in representation theoretic terms with full proofs and \cite[Chapter 10]{GH} for a concise summary.  We start with a partition
$$n = d_1 + d_2 + \ldots + d_r$$
of ${\rm GL}_n$ into blocks of dimension $d_j \geq 1$.  For each $j$ we choose a divisor $f_j \mid d_j$ and in the case $f_j \geq 2$ a cusp form $u_j$ for the group ${\rm SL}_{f_j}(\Bbb{Z})$ with spectral parameter $\mu_j \in \Bbb{C}^{f_j}$ (satisfying \eqref{Lambda} with $f_j$ in place of $n$).  In addition we choose $r$ imaginary numbers $s_1, \ldots s_r\in i\Bbb{R}$ satisfying $\sum_j d_j s_j = 0$. We call the set of such Eisenstein series $E(d, f)$; it consists of all (discrete) choices of such $u_1, \ldots, u_r$ (for $f_j \geq 2$) and such numbers $s_1, \ldots, s_r$ (that vary continuously).  The case of cusp forms corresponds to $r = 1$, $d_1 = f_1 = n$.  In representation theoretic terms, if $\pi$ is a cuspidal automorphic representation on ${\rm GL}_{f}(\Bbb{A})$ corresponding to the cusp form $u$, this corresponds to  the Speh representation ${\rm Speh}(\pi,d/f)$ of ${\rm GL}_d(\Bbb{A})$ which is the unique irreducible subrepresentation of 
\begin{equation}\label{speh}
	{\rm Ind}_{P_{d/f}(\Bbb{A})}^{{\rm GL}_d(\Bbb{A})} \big(\vert \cdot\vert_{\Bbb{A}}^{\frac{d/f-1}{2}}\pi \otimes\ldots \otimes \vert \cdot \vert_{\Bbb{A}}^{-\frac{d/f-1}{2}}\pi\big) 
\end{equation}
where $P_{d/f}$ is the standard parabolic subgroup associated to the partition $f+\ldots +f= d$. \\
 
\textbf{Example:}  Let $n=3$. If $r=1$ and $d_1 = 3$, we have either $f_1 = 3$ in which case we get cusp forms for ${\rm SL}_3(\Bbb{Z})$, or $f_1 = 1$ in which case we get the constant function. If $r= 2$ with $d_1 = 2$ and $d_2 = 1$, we have either $f_1 = 2$ in which case we get maximal Eisenstein series with a cusp form $u_1$ for ${\rm SL}_2(\Bbb{Z})$, or $f_1 = 1$ in which case we get the maximal degenerate Eisenstein series (Epstein zeta function). If $r=3$ with $d_j = f_j = 1$, we get minimal Eisenstein series.\\

From \eqref{speh} we see that the spectral parameter of an element $\varpi \in E(d, f)$ is
\begin{gather*}
  \Big(\underbrace{\mu_1  + s_1 + \textstyle\frac{1}{2} (\frac{d_1}{f_1} - 1)}_{\in \Bbb{C}^{f_1}}, \underbrace{\mu_1  + s_1 + \textstyle\frac{1}{2} (\frac{d_1}{f_1} - 3)}_{\in \Bbb{C}^{f_1}}, \ldots,  \underbrace{\mu_1  + s_1 + \textstyle\frac{1}{2} (1- \frac{d_1}{f_1} )}_{\in \Bbb{C}^{f_1}},\phantom{+++}\\
  \ldots,\\
  \phantom{+++}\underbrace{\mu_r  + s_r + \textstyle\frac{1}{2} (\frac{d_r}{f_r} - 1)}_{\in \Bbb{C}^{f_r}}, \underbrace{\mu_r  + s_r + \textstyle\frac{1}{2} (\frac{d_r}{f_r} - 3)}_{\in \Bbb{C}^{f_r}}, \ldots,  \underbrace{\mu_r  + s_r + \textstyle\frac{1}{2} (1- \frac{d_r}{f_r} )}_{\in \Bbb{C}^{f_r}},\Big)
\end{gather*}
with $\mu_j = 0$ if $f_j = 1$. For each cusp form $u_j$ we use the Jacquet-Shalika bounds to bound  $\| \Re \mu_j \| \leq 1/2$, so that  for $\varpi \in E(d, f)$  we have
\begin{equation}\label{Re}
 \|\Re \mu_{\varpi}\| \leq\max_{1 \leq j \leq r} \Big( \frac{1}{2}\Big( \frac{d_j}{f_j} - 1\Big) + \delta_{f_j \geq 2} \frac{1}{2} \Big). 
\end{equation}

Next, a cusp form in each block of dimension $d_j$ has at most $f_j$ different entries in its spectral parameter, so out of the differences $\mu_{j, i} - \mu_{j, k}$ with $1 \leq i < k \leq d_j$ at least 
\begin{equation*}
f_j \cdot \frac{1}{2} \frac{d_j}{f_j} \Big( \frac{d_j}{f_j} - 1\Big) =  \frac{d_j}{2} \Big( \frac{d_j}{f_j} - 1\Big) 
\end{equation*}
unordered pairs coincide. Using Lemma \ref{lem1}, we conclude that
\begin{equation}\label{supnormEis}
 \int_{  B(\mu) \cap E(d, f)} |\varpi(z) |^2 d\varpi \ll_z (1 + \| \mu \|)^{\frac{1}{2}n(n-1) -\frac{1}{2} \sum_{j=1}^r d_j(d_j/f_j - 1)}
 \end{equation}
 for $\mu \in i\mathfrak{a}^{\ast} $  and $z \in G$. 
This is a reflection of the fact that degenerate Eisenstein series lie on many Weyl chamber walls and therefore have a (somewhat) smaller sup-norm. 
The bound is in general probably far from optimal; see \cite{Bl} for sup-norm bounds for Eisenstein series in a very special case.  
Using \cite[Theorem 1]{JK2} in combination with a local (upper bound) Weyl law  \cite{Mu} for each of the blocks we can do better on average over $z$, namely
\begin{equation}\label{sup-omega}
\int_{\Omega}  \int_{  B(\mu) \cap E(d, f)} |\varpi(z) |^2 d\varpi  \, dz\ll_{\Omega, \varepsilon} (1+ \| \mu \|)^{\frac{1}{2}\sum_jf_j(f_j - 1) +\varepsilon}
\end{equation}
for each compact $\Omega \subseteq {\rm SL}_n(\Bbb{R})$.  

Finally, the set $\{\varpi \in E(d, f) :\| \mu_{\varpi}\| \leq R\}$ can be covered by 
\begin{equation}\label{balls}
\ll R^{\sum_{j=1}^r (f_j-1) + (r-1)} = R^{-1 + \sum_{j=1}^r f_j}
\end{equation}
balls $B(\mu)$. 

\section{Spherical transforms}\label{spher}
Our  goal is to understand the spherical transforms $\tilde{\chi}_T$  and $\tilde{\psi}_{\delta}$.  
By \eqref{trivphi} we have 
\begin{equation}\label{triv}
\tilde{\psi}_{\delta}(\rho) = 1, \quad \tilde{\chi}_{T}(\rho) = \int_{\| g \| \leq T} dg = c_n \text{vol}(\Gamma \backslash G) T^{n(n-1)}
\end{equation}
where the constant $c_n$, defined in \eqref{en}, is computed in \cite[Appendix 1]{DRS}. 
\begin{lemma}\label{lem2} For $0 < \delta< 1$, $A > 0$  and $\Re \mu \ll 1$ we have
$$\tilde{\psi}_{\delta}(\mu) \ll_A (1+  \delta \| \mu \|)^{-A}.$$
 \end{lemma}

\textbf{Proof.} By \eqref{abel} we have $$\tilde{\psi}_{\delta}(\mu) = \int_A a^{\rho} \int_N \psi_{\delta}(an)\, dn \, a^{-\mu} da$$
where we recall that  $\psi_{\delta}$  is $L^1$ supported in a ball of radius $\delta$ about the identity of the $(\frac{1}{2}n(n+1) - 1)$--dimensional space $AN$, so $\| \psi_{\delta} \|_{\infty} \ll \delta^{1-n(n+1)/2}$. The $N$-integral vanishes unless $n \ll \delta$ and $a - \text{id} \ll \delta$. This gives immediately the trivial  bound $\tilde{\psi}_{\delta}(\mu) \ll 1$ for $\Re \mu \ll 1$. On the other hand, we can use an $(n-1)$-dimensional  local coordinate system about the identity in $A$, so that the $A$-integral looks like
$$\int_{\Bbb{R}_{>0}^{n-1}} \Psi_{\delta}(y_1, \ldots, y_{n-1}) y_1^{\mu_1} \cdots y_{n-1}^{\mu_{n-1}} (y_1 \cdots y_{n-1})^{-\mu_n} dy$$
where $\Psi_{\delta}$ is supported in $y_j  = 1 + O(\delta)$ and $ \mathcal{D} \Psi_{\delta}(y) \ll  \delta^{-(n-1) - k}$ for any differential operator of degree $k$ with constant coefficients. Integrating by parts $k$ times with respect to $y_j$   we obtain the bound $\tilde{\psi}_{\delta}(\mu) \ll (\delta |\mu_j-\mu_n|)^{-k}$. Choosing a different local coordinate system, we can replace $n$ with any other index.  This completes the proof.\\

For $a = (a_1 \ldots, a_n) \in \Bbb{R}^n$ we define
$$G(a) = \prod_{j=1}^n (1 + (\max_i a_i) - a_j )^{-1/2} +  \prod_{j=1}^n (1 +| (\min_i a_i) - a_j |)^{-1/2}.$$

\begin{lemma}\label{lem3} Let $T \geq 1$, $\kappa, \varepsilon > 0$ and $\mu$ satisfying \eqref{Lambda}.

There exists a constant $B \in \Bbb{R}$ depending only on $n$ such that
$$\tilde{\chi}_T(\mu) \ll_{\varepsilon} T^{\frac{1}{2}n(n-1) + n\| \Re \mu\| +\varepsilon } (1 + \| \mu \|)^B. \\
$$
If  $ \|\mu \| \ll T^{2-\kappa}$, then  
$$\tilde{\chi}_T(\mu) \ll_{\kappa, \varepsilon} T^{\frac{1}{2}n(n-1) + n\| \Re \mu\| +\varepsilon } \frac{G(\Im \mu)}{(1+ \| \mu \| )^{ \frac{1}{4}n(n-1) + \frac{1}{2}(1 +\| \Re \mu \|)}}. 
$$
\end{lemma} 

\textbf{Proof.} Again we start with \eqref{abel} and compute first $$   \int_N \chi_T(an) dn = \text{meas}\Big( \Big\{\sum_{j= 1}^n a_i^2  + \sum_{1 \leq i < j \leq n} n_{ij}^2 a_i^2  \leq T^2 \mid n_{ij} \in \Bbb{R}\Big\} \Big).$$
We integrate successively  each  $n_{ij}$  at a time using the formula \cite[3.191.1]{GR} 
$$\int_{-\sqrt{Z}}^{\sqrt{Z}} (Z - x^2y^2)^{\alpha} dx = \frac{Z^{\frac{1}{2} + \alpha  } }{y}\frac{\pi^{1/2}\Gamma(1+\alpha)}{\Gamma(\frac{3}{2} + \alpha  )}$$
for $\alpha,   Z, y \geq 0$. In this way we obtain
$$\int_N \chi_T(an) dn  =\gamma_n \delta_{\| a \|_2 \leq T} \frac{(T^2 - \| a \|_2^2)^{n(n-1)/4}}{a_1^{n-1} a_2^{n-2} \cdots a_{n-1}}, \quad  \gamma_n =  \frac{\pi^{n(n-1)/4}}{\Gamma(1 + \frac{1}{4}n(n-1))},$$ 
so that
$$\tilde{\chi}_T(\mu) = \gamma_n  T^{\frac{1}{2}n(n-1)} \int_{\| a \|_2 \leq T} \Big(1 - \frac{\| a \|_2}{T}\Big)^{n(n-1)/4}   \prod_j a_j^{-\mu_j - \frac{n-1}{2}}  da = \gamma_n  \int_{A}f_{\mu}\Big(\frac{a}{T}\Big)  da$$
with
$$f_{\mu}(a) = \delta_{\| a\| \leq 1}(1 - \| a \|_2^2)^{\frac{1}{4}n(n-1)}  \prod_j a_j^{-\mu_j - \frac{n-1}{2}}$$
(using that $\sum_j \mu_j = 0$). Following \cite{DRS}, it is most convenient to estimate this integral asymptotically by passing to the torus in ${\rm GL}_n^+$. Define
$$F_{\mu}(s) = \gamma_n\int_{\Bbb{R}_{>0}^n}f_{\mu}(y)( \det y)^s  \frac{dy_1}{y_1} \cdots \frac{dy_n}{y_n} .$$
Then by Mellin inversion we have 
$$\tilde{\chi}_T(\mu) =   \gamma_n  \int_{(c)} F_{\mu}(s) T^{ns} \frac{ds}{2\pi i} $$
for some sufficiently large $c > 0$. We compute
$$F_{\mu}(s) =   \gamma_n \frac{\Gamma(1 + \frac{n(n-1)}{4})}{2^n\Gamma(\frac{n}{2}s + 1) }  \prod_{j=1}^n \Gamma\Big( \frac{s}{2} - \frac{\mu_j}{2} - \frac{n-1}{4}\Big)$$
using again \cite[3.191.1]{GR} and the fact that $\sum_j \mu_j  = 0$ (this condition on $\mu$ will be used frequently in following arguments). We conclude that
$$\tilde{\chi}_T(\mu) =  \frac{\pi^{n(n-1)/4}}{2^n} T^{n(n-1)/2} \int_{(c)} \frac{T^{ns}}{\Gamma (\frac{n}{2} s + \frac{n(n-1)}{4}+ 1)} \prod_{j=1}^n  \Gamma\Big(\frac{s - \mu_j }{2}\Big) \frac{ds}{2\pi i} $$
for sufficiently large $c > 0$. Estimating this integral is an elaborate exercise in Stirling's formula. Let us write $\mu_j = m_j + i\tau_j$, $s = \sigma + it$ and assume without loss of generality $\tau_1 \geq \tau_2 \geq \ldots \geq \tau_n$. Since $\sum \tau_j = 0$, we have $\tau_1 \asymp - \tau_n \asymp \| \tau\|$. The exponential behavior of the integrand is given by
$$\exp\Big( - \frac{\pi}{4} \sum_{j=1}^n |t-\tau_j| + \frac{\pi n}{4} |t|\Big) \leq \begin{cases}  \exp( - \frac{\pi}{2} \min( |  \tau_1-t|, | t - \tau_n|)), &  \tau_n \leq t \leq \tau_1,\\ 1, & \text{else.}\end{cases}$$
In particular, there is no exponential increase, but exponential decrease for $t \in [\tau_n, \tau_1]$. The polynomial behavior of the gamma quotient is
$$  \ll_{\sigma} (1 + |t|)^{- \frac{1}{2}n\sigma - \frac{1}{4}n(n-1) - \frac{1}{2}} \prod_{j=1}^n (1 + |t - \tau_j|)^{(\sigma - m_j- 1)/2} $$
away from poles. For $|t| \geq 1+ 2 \| \tau \|$ this is $\ll_{\sigma} |t|^{-\frac{1}{4}(n^2 + n +2)}$, in particular, the integrand is absolutely 
 integrable on every vertical line (not crossing poles). Moreover, for $\sigma \leq -n$, $t \ll \| \tau \|$, the gamma quotient is very coarsely bounded by $( 1 + \| \tau \|)^{\frac{1}{2} n |\sigma|}$. In particular, for $\| \mu \| \leq T^{2-\kappa}$, the integral over a line sufficiently far to the left becomes negligible. We use these considerations in the following estimations. 
 
 The right-most pole  appears\footnote{Note that for $\mu = \rho$ the residue of the right-most pole at $s = (n-1)/2$ yields exactly the asymptotic \cite[(A1.15)]{DRS} as it should} at $s =   \| \Re \mu\|. $ Shifting the contour to $c = \| \Re\mu \| + \varepsilon$ and estimating trivially, we obtain immediately the first part of the lemma.
 
Suppose from now on $\|\mu \| \ll T^{2-\kappa}$. We the shift the contour to the far left and pick up the residues. For notational simplicity let us first assume that the $\mu_j$ are pairwise distinct. The general case follows by a straightforward limit procedure. The residues in $\Re s \geq - K$ are given by
$$2 \sum_{j=1}^n \sum_{0 \leq k \leq (K + \Re \mu_j)/2}  \frac{(-1)^k}{k!} \frac{T^{n \mu_j - 2nk}}{\Gamma(\frac{n}{2}\mu_j + \frac{n(n-1)}{4} +1 - nk)} \prod_{i \not=j} \Gamma\Big(\frac{\mu_j - \mu_i}{2} - k\Big).$$
By the same computation as above, each summand is
\begin{equation}\label{summand}
 \ll \frac{T^{nm_j - 2nk}}{(1 + |\tau_j|)^{\frac{n}{2}m_j + \frac{1}{4}n(n-1) + \frac{1}{2} - nk} } \exp\Big(- \frac{\pi}{2}\min(|\tau_1 - \tau_j|, |\tau_j - \tau_n|)\Big)\prod_{i = 1}^n (1 + |\tau_j - \tau_i|)^{\frac{1}{2}(m_j-m_i - 1 - 2k)}.
 \end{equation}
For $\|\tau \| \ll T^{2-\kappa}$ this is decreasing in $k$, so it suffices to consider the term $k=0$. The expression is also increasing in $m_j$, so we can and will assume $m_j = \| m \|$. 
From the exponential term we can assume that $|\tau_j| \asymp \| \tau \|$ in which case
$$\frac{ \prod_{i = 1}^n (1 + |\tau_j - \tau_i|)^{\frac{1}{2}(m_j-m_i )}}{(1 + |\tau_j|)^{\frac{n}{2}m_j }} \ll  \frac{1}{(1 + |\tau_j|)^{\frac{1}{2}m_j}}  \frac{   (1 + |\tau_j |)^{\frac{1}{2}(n-1)m_j- \frac{1}{2}\sum m_i }}{(1 + |\tau_j|)^{\frac{n-1}{2}m_j }} = \frac{1}{(1 + |\tau_j|)^{\frac{1}{2}m_j}} .$$
We therefore bound  \eqref{summand} by
\begin{equation*} 
 \ll \frac{T^{n\| m \|  }}{(1 + |\tau_j|)^{   \frac{1}{4}n(n-1) + \frac{1}{2}(1 + \| m \| ) } } \exp\Big(- \frac{\pi}{2}\min(|\tau_1 - \tau_j|, |\tau_j - \tau_n|)\Big)\prod_{i = 1}^n (1 + |\tau_j - \tau_i|)^{-1/2}
 \end{equation*}
 which gives the desired bound. This completes the proof. \\

Combining \eqref{conv}, Lemma \ref{lem2} and Lemma \ref{lem3}, we conclude
\begin{cor}\label{cor1} Define $\chi_{T, \delta}$ as in \eqref{convdef}. For $n \geq 2$, $  \varepsilon, \kappa, A > 0$, $$T^{-2+\kappa} \ll  \delta < 1\leq T$$ and $\mu$ satisfying  \eqref{Lambda}  we have 
\begin{equation*} 
\begin{split}
\tilde{\chi}_{T, \delta}(\mu) \ll_{\varepsilon, n}   T^{\frac{1}{2}n(n-1) + n\| \Re \mu\| +\varepsilon } \frac{G(\Im \mu)(1+  \delta \| \mu \|)^{-A}}{(1+ \| \mu \| )^{\frac{1}{4} n(n-1) + \frac{1}{2}(1 + \|\Re \mu \|) }} . \end{split}
\end{equation*}
 \end{cor}
 
We remark that for $\| \Re \mu  \| = 0$ the decay in $\mu$ is sharp and the analysis in Lemma \ref{lem3} could be turned into an asymptotic formula. \Valentin{added}

\section{Lattice point count}
For $0 < \delta < 1 \leq T$ we 
 define 
$$\mathcal{N}_{n, \delta}(T; z, w) :=  \sum_{\gamma \in \Gamma} \chi_{T, \delta}(z^{-1}\gamma w).$$
From \eqref{ineq} and \eqref{pretrace} we conclude
$$\mathcal{N}_n(T; z, w)  \leq \mathcal{N}_{n, \delta}(T(1 + \delta) ; z, w) = \int \tilde{\chi}_{T(1+\delta), \delta}(\mu_{\varpi}) \varpi(z)\overline{\varpi(w)} d\varpi.$$
From the right hand side we extract the $L^2$-normalized constant function corresponding to $\mu_{\varpi} = \rho$. By \eqref{triv} this contributes
$$c_n(T(1+\delta))^{n(n-1)} = c_nT^{n(n-1)} + O(T^{n(n-1)}\delta).$$
Similarly we obtain a lower bound and hence conclude the basic asymptotic
$$\mathcal{N}_n(T;z, w) = c_nT^{n(n-1)}  + O\Big(T^{n(n-1)}\delta+ \int_{\mu_{\varpi} \not= \rho} |\tilde{\chi}_{T(1\pm\delta), \delta}(\mu_{\varpi})|( |  \varpi(z)|^2 + | \varpi(w)|^2)d\varpi\Big).$$
We need to estimate the second term.

\subsection{The general argument for $n \geq 5$}\label{gen}
 We partition the spectrum into parameters $(d, f)$ as in Section \ref{Eisen}, excluding the case $r=1$, $d_1 = n$, $f_1 = 1$, which corresponds to the constant function. We assume that $|\log \delta | \asymp \log T$ and specifically $\delta = T^{-\alpha}$ for $0 < \alpha < 2$ (in order to apply Corollary \ref{cor1}). 

Combining Corollary \ref{cor1} (where we drop the factor $G(\mu)$ for simplicity and also simplify the denominator a bit), \eqref{Re}, \eqref{supnormEis} and \eqref{balls}, we have 
\begin{equation}\label{genbound}
\begin{split}
\int_{E(d, f)}(...) &  \ll_{z, w}  T^{\frac{1}{2}n(n-1) + n  \max_{j}  ( \frac{1}{2} ( \frac{d_j}{f_j} - 1 ) + \delta_{f_j \geq 2} \frac{1}{2}  ) + \varepsilon   }\\
&\quad\quad  \times \Big(1 + \delta ^{ \frac{1}{4}n(n-1) + \frac{1}{2}  - \frac{1}{2}n(n-1) + \frac{1}{2} \sum_j d_j(d_j/f_j-1) + 1 - \sum_j f_j}\Big) \\
& =  T^{\frac{1}{2}n(n-1) + \frac{n}{2}  \max_{j} (  \frac{ d_j}{f_j} - \delta_{f_j = 1}) + \varepsilon   }\Big(1 + \delta ^{ -\frac{1}{4}n(n+1) + \frac{3}{2}     + \frac{1}{2} \sum_j (d_j^2 /f_j-2f_j)}\Big).
\end{split}
\end{equation}
Suppose without loss of generality that $j=1$ is the index at which the maximum in the exponent  is attained. Clearly for all other indices the worst case is $f_j = d_j$, so we are left with analyzing
\begin{equation}\label{analy}
T^{\frac{1}{2}n(n-1) + \frac{n}{2}   (  \frac{ d_1}{f_1} - \delta_{f_1 = 1}) + \varepsilon   }\Big(1 + \delta ^{ -\frac{1}{4}n(n+1) + \frac{3}{2} +   \frac{1}{2} (\frac{d_1^2}{f_1} - 2 f_1 - \sum_{j\geq 2} d_j)}\Big). 
\end{equation}
Before we optimize $f_1$, we treat by hand the case $d_1 = n-1$, $f_1 = 1$ (so that $r=2$, $d_2 = 1$), in which case the preceding display becomes
$$T^{n(n-1) - \frac{n}{2}  + \varepsilon   }\Big(1 + \delta ^{  \frac{1}{4} (n^2 - 5n + 2)}\Big) \ll T^{n(n-1) - \frac{n}{2}  + \varepsilon   }$$
for $n \geq 5$. This error term is certainly acceptable. 

From now on we weaken \eqref{analy} a bit and consider
\begin{equation}\label{weaken}
\begin{split}
& T^{\frac{1}{2}n(n-1) + \frac{n d_1}{2f_1} + \varepsilon   }\Big(1 + \delta ^{ -\frac{1}{4}n(n+1) + \frac{3}{2}    + \frac{1}{2} (\frac{d_1^2}{f_1} - 2 f_1 - \sum_{j\geq 2} d_j)}\Big)\\
 = &  T^{\frac{1}{2}n(n-1) + \frac{n d_1}{2f_1} + \varepsilon   }\Big(1 + \delta ^{ -\frac{1}{4}n(n+1) + \frac{3}{2}    + \frac{1}{2} (\frac{d_1^2}{f_1} - 2 f_1 - n+d_1)}\Big).
\end{split}
\end{equation}
Since the cases $d_1 \in \{ n-1, n\}$, $f_1 = 1$ have been ruled out, we always have $d_1/f_1 \leq \max(n-2, n/2) = n-2$  for $n \ge 5$ and so
$$T^{\frac{1}{2}n(n-1) + \frac{n d_1}{2f_1} + \varepsilon   }\leq T^{n(n-1)  - n/2  + \varepsilon}$$
which is clearly  acceptable.  For the second term we need to analyze
$$\phi(\alpha, n, d_1, f_1) = \frac{1}{2}n(n-1) + \frac{n d_1}{2f_1}  - \alpha\Big(-\frac{1}{4}n(n+1) + \frac{3}{2}    + \frac{1}{2} \Big(\frac{d_1^2}{f_1} - 2 f_1 - n+d_1\Big)\Big).$$
We compute
$$\frac{\partial }{\partial d} 
\phi(\alpha, n, d, f)  = \frac{n - \alpha(2d +f)}{2f}$$
with a unique zero at $d_0 = d_0(f)  = (n- \alpha f)/(2\alpha)$ which is a local maximum.  If $n/f \leq \alpha$, then $d_0 \leq  0$, so on the interval $[1, n]$ the function $d \mapsto \phi(\alpha, n, d, f) $ has its maximum at $d = 1$. If $n/f > \alpha$ and $\alpha \geq 1/2$, then $d_0 < n$, and so the maximum lies at $d = d_0$. 

Next we compute 
$$\frac{\partial }{\partial f} 
\phi(\alpha, n, d, f)  = \frac{\alpha d^2 + 2\alpha f^2 - d n}{2f^2} .$$
If $n/d \leq \alpha$, this is always non-negative, so $f \mapsto \phi(\alpha, n, d, f) $ is increasing in $f$. If $n/d >\alpha$, this has a unique positive zero at  $   ((dn - \alpha d^2)/(2\alpha))^{1/2}$ which is a local minimum. So in either case, on the interval $[1, d]$, the function $f \mapsto \phi(\alpha, n, d, f) $  is maximized at $f = 1$ or $f=d$. 

We conclude that for $1 \leq f \leq d \leq n$, the function $\phi(\alpha, n, d, f) $ becomes globally  maximal at most at the three points
$$(f, d)  \in \{ (1, 1), (1, d_0(1)), (n/(3\alpha), n/(3\alpha))\}$$
where $f_0 = n/(3\alpha)$ is the solution to $d_0(f_0) = f_0$. Substituting, we obtain 
$$\phi(\alpha, n, d, f) \leq \max\Big(\frac{n^2}{4}(\alpha + 2) + \frac{3\alpha n}{4} - \frac{3\alpha}{2},  \frac{n^2}{4}\Big(\alpha + \frac{1}{2\alpha} + 2\Big) + \frac{3 n}{4} (\alpha - 1) - \frac{3}{8}\alpha\Big).$$
Thus our final error term is $T^{\psi(\alpha, n) + \varepsilon}$ for 
$$\psi(\alpha, n) = \max \Big(n(n-1) - \alpha, \frac{n^2}{4}(\alpha + 2) + \frac{3\alpha n}{4} - \frac{3\alpha}{2},  \frac{n^2}{4}\Big(\alpha + \frac{1}{2\alpha} + 2\Big) + \frac{3 n}{4} (\alpha - 1) - \frac{3}{8}\alpha\Big)$$
where we can freely choose $0 < \alpha < 2$. A final exercise in calculus shows that for $n \geq 5$ the best choice is 
$$\alpha_0 =  \begin{cases}5/\sqrt{77}, & n = 5,\\ \frac{n}{2n - 1 - \sqrt{2n^2 - 10 n - 4}} = 1 + \frac{1}{\sqrt{2}} + O\big(\frac{1}{n}\big), & n > 5,\end{cases} $$
(satisfying $0 < \alpha_0 < 2$) giving
\begin{equation}\label{final}
\psi(\alpha_0, n) = n(n-1) -   \begin{cases} 5(9 - \sqrt{77})/4, & n = 5,\\ \alpha_0, & n > 5. \end{cases}
\end{equation}
This completes the proof of Theorem \ref{thm1} for $n\geq 5$. \\

Needless to say that these estimates are (deliberately) a bit lossy and can be slightly improved, certainly on a scale $O(1/n)$. In the following two subsections we tighten all screws to obtain (``best-possible'') Selberg type exponents. 
 
\subsection{The case $n=3$} According to the parametrization in Section \ref{Eisen} we distinguish four cases:\\

1) $r = 1$, $d_1 = f_1 = 3$ (cusp forms) with the subcases of a tempered cusp form (case 1a) and a non-tempered cusp form (case 1b);

2) $r = 2$, $d_1= f_1   = 2$, $d_2 =f_2 = 1$ (maximal Eisenstein series with a ${\rm GL}_2$ cusp form) with the subcases of a tempered cusp form (case 2a) and a non-tempered cusp form (case 2b);

3)  $r = 2$, $d_1  = 2$, $d_2 =f_2 = f_1 = 1$ (Epstein zeta function);

4) $r = 3$, $d_1 = d_2 = d_3 = f_1 = f_2 = f_3 = 1$ (minimal Eisenstein series).\\

We can combine the cases 1a, 2a, 4 which are all tempered. Combining Corollary \ref{cor1} (again dropping $\| \Re \mu\|$ in the denominator) and Lemma \ref{lem1}, for each of them we obtain
\begin{displaymath}
\begin{split}
& \int_{E(d, f)} |\tilde{\chi}_{T(1\pm\delta), \delta}(\mu_{\varpi})| |  \varpi(z)|^2  d\varpi  \\
&\ll T^{3+\varepsilon}  \int_{\substack{\mu_1 + \mu_2 + \mu_3 = 0\\ \mu_j \in i\Bbb{R}}} \frac{G(\Im \mu) (1 + \delta \|\mu \|)^{-A}}{(1+\| \mu \|)^{2}} (1 + |\mu_1 - \mu_2|)(1 + |\mu_1 - \mu_3|) (1 + |\mu_2 - \mu_3|) |d\mu| \\
&\ll T^{3+\varepsilon} \delta^{-2} .
\end{split}
\end{displaymath}
Here we used that both terms of $G(\Im \mu)$ contain the square roots of precisely two factors of $(1 + |\mu_1 - \mu_2|)$, $(1 + |\mu_1 - \mu_3|)$, $ (1 + |\mu_2 - \mu_3|)$, so that 
$$G(\Im \mu)(1 + |\mu_1 - \mu_2|)(1 + |\mu_1 - \mu_3|) (1 + |\mu_2 - \mu_3|) \ll (1 + \| \mu \|)^2$$
and we are left with two integration variables of effective length $1/\delta$. 

In the cases 1b, 2b, 3 we have $\| \Re \mu \| \leq 1/2$ and by unitarity (cf.\ \eqref{Lambda})  the spectral spectral parameters are of the form  $(\beta + it, -\beta + it, -2it)$ with $0 < \beta \leq 1/2$, $t\in \Bbb{R}$, so that in each of these cases we can estimate
\begin{displaymath}
\begin{split}
& \int_{E(d, f)} |\tilde{\chi}_{T(1\pm\delta), \delta}(\mu_{\varpi})| |  \varpi(z)|^2  d\varpi  \\
&\ll T^{9/2+\varepsilon}  \int_{t \in \Bbb{R}} \frac{G((t, t, -2t)) (1 + \delta |t|)^{-A}}{(1+|t|)^{2}} (1 + |t|)^2 d t \ll T^{9/2+\varepsilon} \delta^{-1/2} .
\end{split}
\end{displaymath}

Thus we obtain a total error of
$$T^6\delta + T^{3+\varepsilon} \delta^{-2} +  T^{9/2+\varepsilon} \delta^{-1/2} \ll T^{5+\varepsilon}$$
upon choosing $\delta = 1/T$. 
 
\subsection{The case $n=4$} According to the parameterization in Section \ref{Eisen} we distinguish 10 cases:\\
 
1) $r = 1$, $d_1 = f_1 =  4$ (cusp forms) with the subcases of a tempered cusp form (case 1a), a   cusp form with exactly one pair of non-tempered components (case 1b) and a cusp form with 2 pairs of non-tempered components (case 1c);

2) $r = 1$, $d_1 = 4$, $f_1 = 2$ (Speh representation);

3) $r=2$, $d_1 = f_1 = 3$, $d_2 = f_2 = 1$ (maximal Eisenstein series) with the subcases of a tempered ${\rm GL}(3)$ cusp form (case 3a), and a non-tempered cusp form (case 3b);

4) $r=2$, $d_1 = 3$, $f_1 = d_2 = f_2 = 1$ (Epstein zeta function)  

5) $r = 2$, $d_1 = d_2 = f_1 = f_2 = 2$ with the subcases of two tempered ${\rm GL}(2)$ cusp forms (case 5a), exactly one tempered cusp form (case 5b) and two non-tempered cusp forms (case 5c);

6) $r = 2$, $d_1 = d_2 = f_1 = 2$, $f_2 = 1$ with the  subcases of a tempered ${\rm GL}(2)$ cusp form (case 6a) and a non-tempered cusp form (case 6b);

7) $r = 2$, $d_1 = d_2  = 2$, $f_1 = f_2 = 1$;

8) $r = 3$, $d_1 = f_1 = 2$, $d_2 = d_3 = f_2 = f_3 = 1$  with the subcases of a tempered ${\rm GL}(2)$ cusp form (case 8a) and a non-tempered cusp form (case 8b);

9) $r = 3$, $d_1  = 2$, $d_2 = d_3 = f_1 = f_2 = f_3 = 1$;

10) $r = 4$, $d_j = f_j = 1$ (minimal Eisenstein series).\\

We can combine the tempered cases 1a, 3a, 5a, 6a, 8a, 10 and estimate each of them by
\begin{displaymath}
\begin{split}
 &\ll T^{6+\varepsilon}  \int_{\substack{\mu_1 + \mu_2 + \mu_3 + \mu_4 = 0\\ \mu_j \in i\Bbb{R}}} \frac{G(\Im \mu) (1 + \delta \|\mu \|)^{-A}}{(1+\| \mu \|)^{3.5}} \prod_{1 \leq i < j \leq 4} (1 + |\mu_i - \mu_j|) |d\mu| \ll T^{6+\varepsilon} \delta^{-4} 
\end{split}
\end{displaymath}
since $G(\Im \mu)  \prod (1 + |\mu_i - \mu_j|) \ll (1 + \| \mu \|)^{4.5}$. \\

Next we consider the cases 1b, 3b, 5b, 8b, 9. In each of these we have $\| \Re \mu  \| \leq 1/2 $ and spectral parameters of the form $(  \beta + it_1 - \beta +i t_1, -i(t_1 - t_2), -i(t_1 + t_2))$ with $0 \leq \beta \leq 1/2$, $t_1, t_2 \in \Bbb{R}$. This is the only case where we need the extra $\| \Re \mu \| $ in the exponent of the denominator in Corollary \ref{cor1}. 
Here we estimate
\begin{displaymath}
\begin{split}
 &\ll \sup_{0 < \beta \leq 1/2} T^{6 + 4\beta+\varepsilon}  \int_{t_1, t_2 \in \Bbb{R}} \frac{G(t_1, t_1, -t_1 + t_2, - t_1 - t_2) (1 + \delta \| t \|)^{-A}}{(1+\| t \|)^{\frac{7}{2}  +\frac{\beta}{4} }}\\
 &\quad\quad\quad  \times (1 + |2t_1 - t_2|)^2(1 + |2t_1 + t_2|)^2 (1 + |t_2|) dt \ll \sup_{0 < \beta \leq 1/2}  T^{6 + 4\beta +\varepsilon} \delta^{-\frac{5}{2} + \frac{\beta}{4}} .
\end{split}
\end{displaymath}
Here we used that  both terms of $G(...)$ contain the square roots of at  two factors of five linear forms $(1 + |2t_1 - t_2|)^2(1 + |2t_1 + t_2|)^2 (1 + |t_2|)$, and hence $$G(t_1, t_1, -t_1 + t_2, - t_1 - t_2) (1 + |2t_1 - t_2|)^2(1 + |2t_1 + t_2|)^2 (1 + |t_2|)  \ll (1 + \| t \|)^4.$$
Since $\delta \gg T^{-2}$, the worst case is clearly $\beta = 1/2$. \\

Slightly simpler are the cases 1c, 2, 5c, 6b, 7 where again $\| \Re \mu  \| \leq 1/2 $ and the spectral parameters are even more degenerate of the form $(\beta_1 + it, - \beta_1 + it, \beta_2 - it, - \beta_2 - it)$ with $0 <\beta_1, \beta_2 \leq 1/2$, $t \in \Bbb{R}$. Here we have the bound
\begin{displaymath}
\begin{split}
 &\ll T^{8+\varepsilon}  \int_{t \in \Bbb{R}} \frac{G(t, t, -t, -t) (1 + \delta | t |)^{-A}}{(1+| t |)^{3.5}} (1 + |t|)^4 \, dt \ll T^{8 + \varepsilon}  \delta^{-1/2}.
\end{split}
\end{displaymath}\\

It remains to treat case 4 where $\| \Re \mu \| = 1$, and the spectral parameter is of the form $(1 + it, it, -1 + it, -3it)$, so that we obtain the bound
\begin{displaymath}
\begin{split}
 &\ll T^{10+\varepsilon}  \int_{t \in \Bbb{R}} \frac{G(t, t, t, -3t) (1 + \delta | t |)^{-A}}{(1+| t |)^{3.5}} (1 +  |t|)^3  \, dt \ll T^{10 + \varepsilon} . 
\end{split}
\end{displaymath}

Combining the previous bounds, we obtain a total error of
$$\ll T^{12} \delta + T^{6+\varepsilon} \delta^{-4} + T^{8 + \varepsilon} \delta^{-9/4} + T^{10+ \varepsilon} \ll T^{12 - 6/5 + \varepsilon}$$
upon choosing $\delta = T^{-6/5}$. Note that this estimate is very tight and the estimate in the cases 1b, 3b, 5b, 8b, 9 just suffices. \\


\subsection{Proof of Theorem \ref{thm2}} The strategy is the same as in Subsection \ref{gen} except that we replace \eqref{supnormEis} with \eqref{sup-omega}, so that in place of \eqref{genbound} we have to estimate
\begin{equation*}
\begin{split}
  T^{\frac{1}{2}n(n-1) + n  \max_{j}  ( \frac{1}{2} ( \frac{d_j}{f_j} - 1 ) + \delta_{f_j \geq 2} \frac{1}{2}  ) + \varepsilon   }  \Big(1 + \delta ^{ \frac{1}{4}n(n-1) + \frac{1}{2}  -  \frac{1}{2} \sum_jf_j(f_j-1) + 1 - \sum_j f_j}\Big). \\
\end{split}
\end{equation*}
The optimization procedure is again somewhat tedious. We assume that the maximum is attained at $j=1$. Then clearly for all indices $j \geq 2$ the worst case is $f_j = d_j$. For fixed $d_1$, the remaining sum $\sum_{j \geq 2} d_j = n-d_1$ is fixed, so that 
$$\frac{1}{2} \sum_{j =2}^r d_j(d_j+1) $$
becomes maximal if $r = 2$ (the degenerate case $d_1 = n$ would formally correspond to $r=1$). Thus we bound the previous expression by
$$ T^{\frac{1}{2}n(n-1) + \frac{n}{2}     ( \frac{d_1}{f_1}  - \delta_{f_1 = 1}  ) + \varepsilon   }  \Big(1 + \delta ^{ \frac{1}{4}n(n-1) + \frac{3}{2}  -  \frac{1}{2}f_1(f_1 + 1) - \frac{1}{2}(n - d_1)(n - d_1 + 1)   }\Big) .$$
We weaken $f_1(f_1 + 1)$ to $f_1^2 + d_1$, and consider the ``exponent'' function
\begin{displaymath}
\begin{split}
\tilde{\phi}(\alpha, n, d, f) =& \frac{1}{2}n(n-1) + \frac{n}{2}     \Big( \frac{d_1}{f_1}  - \delta_{f_1 = 1}  \Big)  \\
&+ \max\Big[0, - \alpha\Big(\frac{1}{4}n(n-1) + \frac{3}{2}  -  \frac{1}{2}(f_1^2 + d_1)  - \frac{1}{2}(n - d_1)(n - d_1 + 1) \Big)\Big]
\end{split}
\end{displaymath}
where as before $\delta = T^{-\alpha}$ with $0  < \alpha < 2$. It is easy to see that the function $ f \mapsto nd/(2f) - f^2/2$ has its maximum at the boundary, so the worst case 
options for $f_1$ are $f_1 \in \{1, 2, d_1\}$. We have
\begin{displaymath}
\begin{split}
\tilde{\phi}(\alpha, n, d, 1) =& \frac{n(n+d -2) }{2}   + \max\Big[0, - \alpha\Big(\frac{n(n-1)}{4} + 1 -  \frac{d}{2}  - \frac{(n - d)(n - d + 1)}{2} \Big)\Big],\\
\tilde{\phi}(\alpha, n, d, 2) =& \frac{n(n + \frac{1}{2}d -1)}{2}       + \max\Big[0, - \alpha\Big(\frac{n(n-1) }{4}- \frac{1}{2}  -  \frac{d}{2}  - \frac{(n - d_1)(n - d + 1) }{2}\Big)\Big],\\
\tilde{\phi}(\alpha, n, d, d) =& \frac{n^2}{2} + \max\Big[0, - \alpha\Big(\frac{n(n-1) }{4}+ \frac{3}{2}  -  \frac{d^2 + d}{2}  - \frac{(n - d)(n - d + 1)}{2} \Big)\Big]
\end{split}
\end{displaymath}
(the last case if $d > 1$, while   $d=1$ is implicit in the first case). All three functions are non-concave  as functions of $d$, so the maximum can be attained only at the boundary and it suffices to consider
\begin{displaymath}
\begin{split}
\tilde{\phi}(\alpha, n, 1, 1) &= \max\Big[ \frac{n(n-1)}{2}, \frac{2 + \alpha}{4} (n^2 - n) - \frac{\alpha}{2}\Big],\\
  \tilde{\phi}(\alpha, n, n-1, 1) &=  \max\Big[n^2 - \frac{3}{2}n, \frac{4 - \alpha}{4}n^2 - \frac{6 - \alpha}{4} n+ \frac{\alpha}{2} \Big], \\ 
  \tilde{\phi}(\alpha, n, 2, 2) &= \max\Big[\frac{2n^2 - 1}{4}, \frac{2 + \alpha}{4} n^2 - \frac{5\alpha}{4}(n- 2) \Big], \\ 
  \tilde{\phi}(\alpha, n, n, 2) & = \max\Big[\frac{3n^2 - 2n}{4}, \frac{3 - \alpha}{4} n^2 + \frac{\alpha}{4} (n - 2) + \frac{3\alpha}{2} \Big],\\  
  \tilde{\phi}(\alpha, n, n, n) &= \max\Big[\frac{n^2}{2},  \frac{2 + \alpha}{4} n^2 + \frac{3\alpha}{4}(n- 2)\Big].
\end{split}
\end{displaymath}
For $n$ sufficiently large and $1 < \alpha < 2$, the maximum of these values is
$$\max\Big[n^2 - \frac{3}{2}n, \frac{2 + \alpha}{4} n^2 + \frac{3\alpha}{4}(n- 2)\Big],$$
so that our final error term becomes
$$(T^{n(n-1) - \alpha} + T^{n^2 - 3n/2} + T^{\frac{2 + \alpha}{4} n^2 + \frac{3\alpha}{4}(n- 2)})T^{\varepsilon}.$$
The optimal choice for $\alpha$ is
$$\alpha = \frac{2(n^2 - 2n)}{n^2 + 3n - 2} = 2 + O\Big(\frac{1}{n}\Big)$$
as desired.


\begin{thebibliography}{GNY}

\bibitem[Bl]{Bl} V. Blomer, \emph{Epstein zeta-functions, subconvexity, and the purity conjecture},
J. Inst. Math. Jussieu \textbf{19} (2020), 581-596

\bibitem[BM]{BM} V. Blomer, P. Maga, \emph{Subconvexity for sup-norms of cusp forms on ${\rm PGL}(n)$},
Selecta Math. \textbf{22} (2016), 1269-1287

\bibitem[DRS]{DRS} W. Duke, Z. Rudnick, P. Sarnak, \emph{Density of integer points on affine homogeneous varieties}, Duke Math. J. \textbf{71} (1993),  143-179

\bibitem[GH]{GH} J. Getz, H. Hahn, \emph{An introduction to automorphic representations} (2019)



\bibitem[GNY]{GNY} A. Gorodnik, A. Nevo, G. Yehoshua, \emph{Counting lattice points in norm balls on higher rank simple Lie groups}, Math. Res. Lett. \textbf{24} (2017),  1285-1306.

\bibitem[GR]{GR} I. S.  Gradshteyn, I. M. Ryzhik, \emph{Table of integrals, series, and products}, 7th ed.,  Academic Press 2007

\bibitem[He]{He} S. Helgason, \emph{Groups and geometric analysis. Integral geometry, invariant differential operators, and spherical functions},  Mathematical Surveys and Monographs  \textbf{83}. American Mathematical Society, Providence, RI, 2000.

\bibitem[Iw]{Iw} H. Iwaniec, \emph{Spectral methods of automorphic forms}, Grad. Stud. Math. \textbf{53} (2002), AMS

\bibitem[IK]{IK}  H. Iwaniec, E. Kowalski, \emph{Analytic Number Theory},  AMS Colloquium Publications \textbf{53} (2004)

\bibitem[JK1]{JK2} S. Jana, A. Kamber, \emph{On the local $L^2$-bound of the Eisenstein series}, {\tt arXiv:2210.16291}

\bibitem[JK2]{JK} S. Jana, A. Kamber, \emph{Optimal diophantine exponents for ${\rm SL}(n)$}, {\tt arXiv:2211.05106}

\bibitem[LP]{LP} P. Lax, R. Phillips, \emph{The asymptotic distribution of lattice points in Euclidean and non-Euclidean spaces}, 
J. Funct. Anal. \textbf{46} (1982),  280-350.

\bibitem[MW1]{MW1}  C. M{\oe}glin, J.-L. Waldspurger, \emph{Le spectre r\'esiduel de ${\rm GL}(n)$}, Ann. Sci \'Ecole Norm. Sup. \textbf{22} (1989), 605-674


\bibitem[MW2]{MW2}  C. M{\oe}glin, J.-L. Waldspurger, \emph{Spectral decomposition and Eisenstein series}, Cambridge Tracts in Mathematics \textbf{113} (1995)

\bibitem[Mu]{Mu} W. M\"uller, \emph{Weyl's law for the cuspidal spectrum of ${\rm SL}_n$}, Ann. of Math.  \textbf{165} (2007),  275-333

\bibitem[PR]{PR} Y. Petridis, M. Risager, \emph{Local average in hyperbolic lattice point counting}, with an appendix by N. Laaksonen, Math. Z. \textbf{285} (2017), 1319-1344


\bibitem[Se]{Sel} A. Selberg, \emph{Harmonic analysis and discontinuous groups in weakly symmetric Riemannian spaces with applications to Dirichlet series}, J. Indian Math. Soc. B \textbf{20} (1956), 47-87


\bibitem[Sz]{Sz}  W. Szerpin\'ski, \"Uber ein Problem aus der analytischen Zahlentheorie, Prace mat.-fiz. \textbf{17}
(1906), 77-118

\end{thebibliography}
\end{document}